\documentclass[onecolumn,12pt]{IEEEtran} 
\usepackage{amssymb}
\usepackage{amsthm,color}
\usepackage{latexsym}
\usepackage{graphicx}
\usepackage{fancyhdr}
\usepackage{bbm}
\usepackage{amsmath}
\usepackage{latexsym,amssymb,amsxtra,mathrsfs,times}

\renewcommand{\vec}[1]{{\bf #1}}

\newtheorem{thm}{Theorem}

\newtheorem{lemma}[thm]{Lemma}

\theoremstyle{definition}
\newtheorem{defin}[thm]{Definition}

\newcommand{\tr}{{\rm Tr}}
\newcommand{\R}{\mathbb R}
\newcommand{\G}{\mathcal{G}(s)}

\newcommand{\F}{{\rm GF}}
\newcommand{\C}{{\mathcal C}}
\newcommand{\CC}{{\mathbb C}}

\newcommand{\D}{{\mathcal D}}
\newcommand{\B}{{\mathcal B}}

\newcommand{\Ea}{{\mathbb{E}\left(A_l(s), \Omega_p\right)}}
\newcommand{\E}{{\mathbb{E}}}

\newcommand{\A}{{A_l(s)}}

\newcommand{\MP}{\mathrm{MP}}

\newcommand{\ud}{\mathrm{d}}

\usepackage{setspace}
\begin{document}




\title{On a question of Babadi and Tarokh}

\author{Jing Xia$^1$ \thanks{1. Fred Hutchinson Cancer Research Center, 1100 Fairview Ave N, Seattle, WA, USA}, Maosheng Xiong$^2$ \thanks{2. Department of Mathematics, Hong Kong University of Science and Technology, Clear Water Bay, Kowloon, Hong Kong}}





\maketitle


\renewcommand{\thefootnote}{}





\begin{abstract}
In a recent remarkable paper \cite{ba1}, Babadi and Tarokh proved the ``randomness'' of sequences arising from binary linear block codes in the sense of spectral distribution, provided that their dual distances are sufficiently large. However, numerical experiments conducted by the authors revealed that Gold sequences which have dual distance 5 also satisfy such randomness property. Hence the interesting question was raised as to whether or not the stringent requirement of large dual distances can be relaxed in the theorem in order to explain the randomness of Gold sequences. This paper improves their result on several fronts and provides an affirmative answer to this question.
\end{abstract}


\section{Introduction}\label{sec-into}

The elegant theory of random matrices, and in particular properties of their spectral distribution, have been studied for a long time but remain a prominent and active research area due to its wide and important applications in many diverse disciplines such as mathematical statistics, theoretical physics, number theory, and more recently in economics \cite{paf} and communication theory \cite{tul}. Most of the random models considered so far are matrices whose entries have i.i.d. structures. In a remarkable paper, Babadi and Tarokh \cite{ba1} considered matrices formed by choosing randomly codewords from some linear block codes with large dual distance and proved that these matrices behave like random matrices with i.i.d. entries, as long as the empirical spectral distribution is concerned. To describe their beautiful result, we need some notation.

Let $\C$ be an $[n,k,d]$ binary linear block code of length $n$, dimension $k$ and minimum Hamming distance $d$ over $\F(2)$. The dual code of $\C$, denoted by $\C^{\bot}$, is an $[n,n-k,d^{\bot}]$ binary linear block code over $\F(2)$ such that all the codewords of $\C^{\bot}$ are orthogonal to those of $\C$ with the inner product defined over $\F(2)^n$. Let $\epsilon: \F(2)^n \to \{-1,1\}^n$ be the component-wise mapping $\epsilon(v_i):=(-1)^{v_i}$, for $\vec{v}=(v_1,v_2,\ldots,v_n) \in \F(2)^n$. For $p<n$, let ${\bf \Phi}_{\C}$ be a $p \times n$ random matrix whose rows are obtained by mapping a uniformly drawn set of size $p$ of the codewords of $\C$ under $\epsilon$. The \emph{Gram matrix} of ${\bf \Phi}_{\C}$ is defined as $\mathcal{G}_{\C}:={\bf \Phi}_{\C}{\bf \Phi}_{\C}^T$, where ${\bf \Phi}_{\C}^T$ is the transpose of ${\bf \Phi}_{\C}$. Let $\{\lambda_1,\lambda_2,\ldots,\lambda_n\}$ be the set of eigenvalues of an $n \times n$ matrix $\vec{A}$. The \emph{spectral measure} of $\vec{A}$ is defined by
\[\mu_{\vec{A}}:=\frac{1}{n} \sum_{i=1}^n \delta_{\lambda_i}, \]
where $\delta_z$ is the Dirac measure. The empirical spectral distribution of $\vec{A}$ is defined as
\[M_{\vec{A}}(z):=\int_{-\infty}^z \mu_{\vec{A}}(\ud z). \]
Babadi and Tarokh proved the following result (\cite[Theorem 2.1]{ba1}):

\emph{Consider a sequence of $[n,k_n,d_n]$ binary linear block codes $\{\C_n\}_{n=1}^{\infty}$. Let ${\bf \Phi}_{\C_n}$ be a $p \times n$ random matrix based on $\C_n$, let $\mathcal{G}_{\C_n}$ denote the Gram matrix of the matrix $\frac{1}{\sqrt{n}}{\bf \Phi}_{\C_n}$, and let $M_{\C_n}(z)$ denote the empirical spectral distribution of $\mathcal{G}_{\C_n}$. Finally, let $r_n$ be the greatest even integer less than or equal to $[(d^{\bot}_n-1)/2]$, and let $r:=\liminf_n r_n$. Then, as $n \to \infty$ with $y:=p/n \in (0,1)$ fixed, we have}
\[\limsup_n |M_{\C_n}(z)-M_{\MP}(z)| \le c(y,r) \left(r^{-1}+r^{-2}\right)\]
\emph{almost surely for all $z$, where $c(y,r)$ is a bounded function of $r$ (which can be given explicitly), and $M_{\MP}(z)$ is the distribution corresponding to the Marchenko-Pastur measure $\mu_{\mathrm{MP}}$ whose density is given by }
\[\frac{\ud \mu_{\MP}}{\ud z}:=\frac{1}{2 \pi z y} \sqrt{(b-z)(z-a)} \, 1_{(a \le z \le b)}\,\,, \]
\emph{here $a=(1-\sqrt{y})^2$ and $b=(1+\sqrt{y})^2$.}

It is well-known that as the dimensions grow to infinity, the empirical spectral distribution of the Gram matrix of real i.i.d. random matrices follows the Marchenko-Pastur law \cite{mar}. With this respect, the above result indicates that the matrix $\frac{1}{\sqrt{n}}{\bf \Phi}_{\C}$ based on the binary linear block code $\C$ is very close to random i.i.d. generated matrices as $n \to \infty$, if the dual distance of the code $\C$ is large enough. Numerical experiments conducted by the authors \cite{ba1} on some low-rate BCH codes confirmed the significant similarity of the empirical distribution to the Marchenko-Pastur law for dimensions (and consequently, dual distances) as small as $n=63$.

However, there is an interesting phenomenon: the authors \cite{ba0} also conducted some numerical experiments on Gold sequences and found convincing similarity of the empirical distributions to the Marchenko-Pastur law as well. This is a little surprising because Gold sequences arise from Gold codes \cite{G} whose dual distances are always 5, which is relatively small. In a more recent interesting paper \cite{ba2}, investigating much further on the topic, the authors proved decisively the ``randomness'' of products of matrices arising from different binary linear block codes under large dual distances. At the end of the paper \cite{ba2} Babadi and Tarokh also conducted numerical experiments and found numerical evidence of randomness on some Gold sequences. Hence they raised the natural question as to relaxing the stringent requirement of large dual distances in the results in order to explain the mysterious randomness of Gold sequences.

The purpose of this paper is to provide an affirmative answer to this questions. While binary linear block codes are most useful in practice, it is worthwhile to consider, at least in theory, linear block codes over a general finite field $\F(q)$ where $q$ is a prime power, especially when it does not require any substantial effort. For this purpose, denote by $\psi: \F(q) \to \CC^*$ the standard additive character given by
\[\psi(z)=\exp\left(\frac{2 \pi \sqrt{-1} \, \tr_{q/l}(z)}{l}\right),\]
here $l$ is any prime number and $q$ is a power of $l$, and $\tr_{q/l}$ denotes the trace mapping from $\F(q)$ to $\F(l)$. When $q=l=2$, then $\psi(z)=(-1)^z$ for $z \in \F(2)$ which was considered before. It is known that $\psi(z)$ is a complex $p$-th root of unity.

Let $\C$ be an $[n,k,d]$ linear block code of length $n$, dimension $k$ and minimum Hamming distance $d$ over $\F(q)$. The dual code of $\C$, denoted by $\C^{\bot}$, is an $[n,n-k,d^{\bot}]$ linear block code over $\F(q)$ such that all the codewords of $\C^{\bot}$ are orthogonal to those of $\C$ with the natural inner product defined over $\F(q)^n$. Let $\epsilon: \F(q)^n \to (\CC^*)^n$ be the component-wise mapping $\epsilon(v_i):=\psi(v_i)$, for $\vec{v}=(v_1,v_2,\ldots,v_n) \in \F(q)^n$. For $p<n$, let ${\bf \Phi}_{\C}$ be a $p \times n$ random matrix whose rows are obtained by mapping a uniformly drawn set of size $p$ of the codewords of $\C$ under $\epsilon$. The \emph{Gram matrix} of the $p \times n$ matrix ${\bf \Phi}_{\C}$ is defined as $\mathcal{G}_{\C}:={\bf \Phi}_{\C}{\bf \Phi}_{\C}^*$, where ${\bf \Phi}_{\C}^*$ is the conjugate transpose of ${\bf \Phi}_{\C}$. We prove

\begin{thm} \label{1:thm} Let $\C$ be an $[n,k,d]$ linear block code over $\F(q)$. Let ${\bf \Phi}_{\C}$ be a $p \times n$ random matrix based on $\C$, let $\mathcal{G_{\C}}$ denote the Gram matrix of $\frac{1}{\sqrt{n}}{\bf \Phi}_{\C}$, and let $M_{\C}(z)$ denote the empirical spectral distribution of $\mathcal{G}_{\C}$. Suppose $n$ is sufficiently large. Then if $d^{\bot} \ge 5$ and for any $y:=p/n \in (0,1)$, we have
\begin{eqnarray} \label{1:eqnthm} \sup_{z \in \R}|M_{\C}(z)-M_{\MP}(z)| \le \frac{800}{\sqrt{y}(1-y)}\, \frac{\log \log n}{\log n} \,\,. \end{eqnarray}
\end{thm}

\subsection{Discussion of the Main Theorem}

Theorem \ref{1:thm} might look a little surprising, compared with the celebrated result by Sidel'nikov \cite{sid}: for any $[n,k,d]$ binary linear block code $\C$ with $d^{\bot} \ge 3$, we have
\[A(z)-\Phi(z)| \le \frac{9}{\sqrt{d^{\bot}}} \]
as $n \to \infty$, where $A(z)$ is the cumulative weight distribution function of the code $\C$ and
\[\Phi(z):=\frac{1}{\sqrt{2 \pi}} \int_{-\infty}^{z}e^{-t^2/2} \ud t\,. \]
Hence the ``randomness'' of the weight distribution of $\C$ is ensured if $d^{\bot}$ is sufficiently large. In Theorem \ref{1:thm}, however, we only require $d^{\bot} \ge 5$.

Gold codes have three distinct non-zero weights which are known \cite{G}. By applying the MacWilliams identity \cite{mac} and by using {\bf Mathematica}, it can be readily verified that the dual distance of Gold codes is always 5, hence Theorem \ref{1:thm} is applicable and confirms that Gold sequences behave like random i.i.d. sequences, in the sense of the spectral distribution.

The condition $d^{\bot} \ge 5$ in Theorem \ref{1:thm} can be slightly improved by assuming that the number of weight 4 codewords in $C^{\bot}$ is relatively small (see Theorem \ref{2:moment} in Section \ref{sec-moment}), and the inequality (\ref{1:eqnthm}) of same kind still holds true, if $800$ replaced by a larger constant on the right hand side of (\ref{1:eqnthm}). On the other hand, however, if $d^{\bot}=3$, then Theorem \ref{1:thm} may not be true: Babadi, Ghassemzadeh and Tarokh (\cite[Theorem 3.1]{ba0}) proved that shortened first-order Reed-Muller (Simplex) codes which have dual distance 3 have substantially different behavior in the sense of the spectral distribution.

The proof of Theorem \ref{1:thm} follows essentially the strategy used by Babadi and Tarokh in \cite{ba1}, but here in the paper some essence of number theory plays more prominent roles in the study. This might become more apparent in Section \ref{sec-moment} when we study the $l$-moment of the spectral measure. We shall prove Theorem \ref{2:moment}, which improves \cite[Lemma 3.3]{ba1} substantially. Equipped with Theorem \ref{2:moment}, in Section \ref{sec-thm1} we will prove Theorem \ref{1:thm} directly. In the proof of Theorem \ref{2:moment}, however, some very complicated issues of combinatorial nature arise which need to be taken care of. To streamline the ideas of the paper, we treat those issues in Section \ref{sec-combin}. 

\section{Estimate of the $l$-th moment}\label{sec-moment}

In this section we study the $l$-th moment of the spectral distribution, similar to \cite[Lemma 3.3]{ba1}. We use slightly different notation, which might be more suited for the problem.

As in Introduction, let $\C$ be an $[n,k,d]$ linear block code over $\F(q)$, and let $\epsilon: \F(q)^n \to (\CC^*)^n$ be the component-wise mapping. Define $\D=\epsilon(\C)$. Let $N:=q^k$ be the cardinality of $\D$ (and $\C$). Let $p<n$. In order to choose randomly $p$ elements from $\D$, we define $\Omega_p$ to be the set of all maps $s: [1,p] \to \D$ endowed with the uniform probability, here $[1,p]$ denotes the set of integers from $1$ to $p$. Hence $\Omega_p$ is a probability space with cardinality $|\Omega_p|=N^p$. For each $s \in \Omega_p$, the $p \times n$ matrix ${\bf \Phi}_s$ corresponding to $s$ is given by
\[{\bf \Phi}_s^T=\left[s(1)^T,s(2)^T, \ldots,s(p)^T \right]_{n \times p}\,\,, \]
here we have written $s(i) \in \D$ as $1 \times n$-row vectors. For any $\vec{u}=(u_1,\ldots,u_n),\vec{v}=(v_1,\ldots,v_n) \in \CC^n$, the (Hermitian) inner product is
\[\langle \vec{u},\vec{v} \rangle:=u_1 \bar{v}_1+\cdots+u_n \bar{v}_n.\]
Let $\mathcal{G}(s)$ be the Gram matrix of $\frac{1}{\sqrt{n}} {\bf \Phi}_s$. This is a $p \times p$ Hermitian matrix with the $(ij)$-th entry given by $\langle s(i),s(j)\rangle /n$. Let $\lambda_1(s), \lambda_2(s), \ldots, \lambda_p(s) \in \R$ be the eigenvalues of $\mathcal{G}(s)$. For any positive integer $l$, define
\[\A:=\frac{1}{p} \sum_{i=1}^p \lambda_i(s)^l=\frac{1}{p}\,\, \tr \left(\G^l\right). \]
The purpose of this section is to compute $\Ea$, the $l$-th moment of the spectral measure. We prove a general result:

\begin{thm} \label{2:moment}
Let $y:=p/n \in (0,1)$. Let $A$ be the number of weight 4 codewords in $\C^{\bot}$. Then for any $2 \le l <\sqrt{p}$, we have
\begin{eqnarray} \label{2:momenteqn}
\Ea=\sum_{i=0}^{l-1}\frac{y^i}{i+1}\binom{l}{i} \binom{l-1}{i}+E_l,
\end{eqnarray}
where $E_l$ is bounded by
\[|E_l| \le \left(4+2\sqrt{\frac{2A}{q-1}+\frac{1}{4}}\right) \frac{l^{l+1}}{n}\,\,, \]
\end{thm}
The rest of this section is devoted to a proof of Theorem \ref{2:moment}.

\subsection{Problem setting up}

We say that $\gamma: [0,l] \to [1,p]$ is a closed path if $\gamma$ is a map with $\gamma(0)=\gamma(l)$. Denote by $\Pi_{l,p}$ the set of all closed paths from $[0,l]$ to $[1,p]$. For each $\gamma \in \Pi_{l,p}$ and $s \in \Omega_p$, define
\[\omega_{\gamma}(s):=\langle s \circ \gamma(0), s \circ \gamma (1)\rangle \langle s \circ \gamma(1), s \circ \gamma (2)\rangle \cdots \langle s \circ \gamma(l-1), s \circ \gamma (l)\rangle. \]
Expanding $\tr\left(\G^l\right)$, it is easy to see that
\[A_l(s)=\frac{1}{pn^l}\sum_{\gamma \in \Pi_{l,p}} \omega_{\gamma}(s). \]
Hence
\[\Ea=\frac{1}{pn^l}\sum_{\gamma \in \Pi_{l,p}} \E\left(\omega_{\gamma}(s), \Omega_p\right). \]
Let $\Sigma_p$ be the group of permutations of the set $[1,p]$. Then $\Sigma_p$ acts on $\Pi_{k,p}$, since $\sigma \circ \gamma \in \Pi_{l,p}$ whenever $\gamma \in \Pi_{l,p}$ and $\sigma \in \Sigma_p$. Let $[\gamma]$ be the equivalent class of $\gamma$, that is,
\[[\gamma]=\{\sigma \circ \gamma: \sigma \in \Sigma_p \}. \]
We may write
\[\Ea=\frac{1}{pn^l}\sum_{\gamma \in \Pi_{l,p}/\Sigma_p} \,\,\, \sum_{\tau \in [\gamma]}\E\left(\omega_{\tau}(s), \Omega_p\right). \]
For any fixed $\sigma \in \Sigma_p$, as $s$ runs over $\Omega_p$, clearly $s \circ \sigma$ also runs over $\Omega_p$, hence
\[\E\left(\omega_{\sigma \circ \gamma}(s), \Omega_p\right)=\E\left(\omega_{\gamma}(s \circ \sigma), \Omega_p\right)=\E\left(\omega_{\gamma}(s), \Omega_p\right). \]
Moreover, let
\[V_{\gamma}:=\gamma \left([0,l]\right) \subset [1,p], \quad v_{\gamma}:=\# V_{\gamma},\]
and define the probability space
\[\Omega(V_{\gamma}):=\{s: V_{\gamma} \to \D\}\]
assigned with the uniform probability. It is clear that $\#[\gamma]=\frac{p!}{(p-v_{\gamma})!}, \#\Omega(V_{\gamma})=N^{v_{\gamma}}$ and
\[\E\left(\omega_{\gamma}(s), \Omega_p\right)=\E\left(\omega_{\gamma}(s), \Omega(V_{\gamma})\right). \]
Summarizing the above considerations, we have
\begin{eqnarray} \label{2:set} \Ea=\frac{1}{pn^l}\sum_{\gamma \in \Pi_{l,p}/\Sigma_p} \frac{p!}{(p-v_{\gamma})!}\,\, \E\left(\omega_{\gamma}(s), \Omega\left(V_{\gamma}\right)\right). \end{eqnarray}

\subsection{Study of $\E\left(\omega_{\gamma}(s), \Omega\left(V_{\gamma}\right)\right)$}

Up to this point everything is essentially the same as in the proof of \cite[Lemma 3.3]{ba1}. The main innovation of the paper is to use number theory to treat $\E\left(\omega_{\gamma}(s), \Omega\left(V_{\gamma}\right)\right)$ in a more careful way.

Let $H=(h_{ij})_{n \times k}$ be a generating matrix of $\C$, that is, each codeword of $\C$ is given by
\begin{eqnarray} \label{2:codec} c(\vec{x}):=H [x_1,\ldots,x_k]^T,  \end{eqnarray}
for some unique $\vec{x}=(x_1,\ldots,x_k) \in \F(q)^k$. Hence each $s(i) \in \D$ corresponds to a unique vector, which we may record as $(s(i)_1,\ldots,s(i)_k) \in \F(q)^k$. From (\ref{2:codec}), the $t$-th entry of $s(i)$ is given by
\[s(i)[t]=\psi\left(\sum_{j=1}^kh_{tj}s(i)_j\right), \]
where $\psi:\F(q) \to \CC^*$ is the standard additive character. So
\[\langle s \circ \gamma (u), s \circ \gamma (u+1) \rangle= \sum_{t=1}^n \psi\left(\sum_{j=1}^kh_{tj}s \circ \gamma(u)_j-\sum_{j=1}^kh_{tj}s \circ \gamma (u+1)_j\right), \]
and hence
\begin{eqnarray*} \omega_{\gamma}(s) &=& \sum_{1 \le t_0,t_1,\ldots,t_{l-1} \le n} \psi\left(\sum_{j=1}^kh_{t_0j} \left\{s \circ \gamma(0)_j-s \circ \gamma (1)_j\right\} \right) \\
& &  \times \psi\left(\sum_{j=1}^kh_{t_1j} \left\{s \circ \gamma(1)_j-s \circ \gamma (2)_j\right\} \right) \cdots \psi\left(\sum_{j=1}^kh_{t_{l-1}j} \left\{s \circ \gamma(l-1)_j-s \circ \gamma (0)_j\right\} \right).
\end{eqnarray*}
Now suppose
\[V_{\gamma}=\{z_a: 1 \le a \le v_{\gamma}\} \subset [1,p],\]
and for each $a$, let $I_a:=\gamma^{-1}(z_a)$. For each $u \in I_a$, we have $\gamma(u)=z_a$ and clearly $[0,l-1]=\bigcup_a I_a$ is a partition. We may collect the term $s(z_a)$ together on the right hand side of $\omega_{\gamma}(s)$ above and rewrite it as
\[\omega_{\gamma}(s)=\sum_{1 \le t_0,t_1,\ldots,t_{l-1} \le n} \prod_{a=1}^{v_{\gamma}} \prod_{j=1}^k\psi \left(s(z_a)_j \sum_{u \in I_a} \{h_{t_u j}-h_{t_{u-1}j}\}\right). \]
Here when $u=0$, we interpret $t_{0-1}:=t_{l-1}$ (we will use this convent multiple times in the paper). Therefore
\begin{eqnarray*} \label{2:ees} \E\left(\omega_{\gamma}(s), \Omega\left(V_{\gamma}\right)\right)=\frac{1}{N^{v_{\gamma}}}\sum_{
\substack{s(z_a)_j \in \F(q)\\
1 \le a \le v_{\gamma}\\
1 \le j \le k}} \omega_{\gamma}(s). \end{eqnarray*}
The orthogonality property
\[\sum_{z \in \F(q)} \psi(zx)=\left\{\begin{array}{lll}
0&:& \mbox{ if } x \in \F(q) \setminus \{0\};\\
q&:& \mbox{ if } x=0,
\end{array}
\right.\]
implies that if for some $a$ and for some $j$ we have
\[\sum_{u \in I_a} \left(h_{t_uj}-h_{t_{u-1}j} \right) \ne 0 \,, \]
then their contribution to $\E\left(\omega_{\gamma}(s), \Omega\left(V_{\gamma}\right)\right)$ is zero. So we conclude that the quantity $\E\left(\omega_{\gamma}(s), \Omega\left(V_{\gamma}\right)\right)$ is the same as $W_{\gamma}$, which is the number of solutions $(t_0,t_1,\ldots,t_{l-1})$ such that $1 \le t_0,t_1,\ldots,t_{l-1} \le n$ and
\[\sum_{u \in I_a} \left(\vec{h}_{t_u}-\vec{h}_{t_{u-1}} \right)=\vec{0}, \quad \forall \, \, 1 \le a \le v_{\gamma},\]
here $\vec{h}_i$ denotes the $i$-th row of the matrix $H$, and
\begin{eqnarray} \label{2:set} \Ea=\frac{1}{pn^l}\sum_{\gamma \in \Pi_{l,p}/\Sigma_p} \frac{p!}{(p-v_{\gamma})!}\, W_{\gamma}. \end{eqnarray}

\subsection{Proof of Theorem \ref{2:moment}}

The combinatorial nature of solving $W_{\gamma}$, while elementary, presents some technical challenge. To streamline the idea of the proof, and for the sake of clarity, we leave the analysis of $W_{\gamma}$ to Section \ref{sec-combin}. Here instead we quote the main results to continue our proof.

In Section \ref{sec-combin} we prove that there is a subset $\Gamma \subset \Pi_{l,p}/\Sigma_p$ with the following property:
\[\left\{\begin{array}{lll} W_{\gamma}= n^{l-v_{\gamma}+1} &: & \mbox{ if } \gamma \in \Gamma; \\
0 \le W_{\gamma} \le C_A \, n^{l-v_{\gamma}} &: & \mbox{ if } \gamma \notin \Gamma, \end{array}
\right. \]
where $C_A$ is given in (\ref{2:ca}). Using this we find that
\begin{eqnarray} \label{2:mainestimate} \Ea=\frac{n}{p}\sum_{\substack{\gamma \in \Pi_{l,p}/\Sigma_p\\
\gamma \in \Gamma}} \frac{p!}{(p-v_{\gamma})! \, n^{v_{\gamma}}}+E_1, \end{eqnarray}
where $E_1$ is bounded by
\[|E_1| \le \frac{C_A}{p}\sum_{\gamma \in \Pi_{l,p}/\Sigma_p} \frac{p!}{(p-v_{\gamma})! \, n^{v_{\gamma}}} \le \frac{C_A}{p}\sum_{v=1}^l \left(\frac{p}{n}\right)^{v} \sum_{\substack{\gamma \in \Pi_{l,p}/\Sigma_p\\
v_{\gamma}=v}} 1\,\, . \]
It is easy to see that
\[\sum_{\substack{\gamma \in \Pi_{l,p}/\Sigma_p\\
v_{\gamma}=v}} 1 < v^l \le l^l, \]
and hence
\[|E_1| \le C_A \, l^{l+1}/n. \]
On the other hand, it is also proved in Section \ref{sec-combin} that
\[\sum_{\substack{\gamma \in \Gamma \subset \Pi_{l,p}/\Sigma_p \\
v_{\gamma}=v}} 1=\frac{1}{v} \binom{l}{v-1} \binom{l-1}{v-1}. \]
Suppose $2 \le l <\sqrt{p}$. For $v \ge 2$, using
\[p^v \ge \frac{p!}{(p-v)!}>p^v \left(1-v/p\right)^{v-1} \ge p^v\left(1-v(v-1)/p\right), \]
in (\ref{2:mainestimate}), we can finally obtain, after some simplifying, the desired result (\ref{2:momenteqn}). This completes the proof of Theorem \ref{2:moment}. \quad $\square$

\section{Proof of Theorem \ref{1:thm}} \label{sec-thm1}

Given Theorem \ref{2:moment}, the proof of Theorem \ref{1:thm} follows essentially arguments in \cite{ba1}, though some of our analysis is more precise.

\subsection{Some lemmas} Fix $y \in (0,1)$, let $\vec{x}$ be a Marchenko-Pastur random variable whose density function is given by
\[\frac{\ud \mu_{\MP}}{\ud z}:=\frac{1}{2 \pi z y} \sqrt{(b-z)(z-a)} \, 1_{(a \le z \le b)}\,\,, \]
\emph{here $a=(1-\sqrt{y})^2$ and $b=(1+\sqrt{y})^2$.}
It is known that the $l$-th moment of $\vec{x}$ is given by
\begin{eqnarray} \label{3:mom} m_{\MP}^{(l)}=\E(\vec{x}^l)=\sum_{i=0}^{l-1}\frac{y^i}{i+1} \binom{l}{i} \binom{l-1}{i}. \end{eqnarray}
Define
\[b_{\MP}^{(l)}:=\E\left((\vec{x}-1)^l\right). \]
Clearly $b_{\MP}^{(0)}=1, b_{\MP}^{(1)}=0$. We first prove

\begin{lemma} \label{3:lem1} For any $l \ge 2$ we have
\begin{eqnarray} \label{3:bmp} \left|b_{\MP}^{(l)}\right| < \frac{l^3 (8e^2)^l y}{8 \pi}\, . \end{eqnarray}
\end{lemma}

\noindent {\bf Proof.} Expanding $\E\left((\vec{x}-1)^l\right)$ and using (\ref{3:mom}) we have
\[b_{\MP}^{(l)}=\sum_{i=1}^{l-1} \frac{y^i}{i+1}\sum_{t=i+1}^l (-1)^{l-1} \binom{l}{t}\binom{t}{i}\binom{t-1}{i}. \]
Elementary estimates on binomial coefficients yield
\[\left|b_{\MP}^{(l)}\right|< \frac{2^l}{2} \sum_{i=1}^{l-1} \frac{y^i l^{2i}}{(i!)^2} < 2^{l-1} (y l^2) \sum_{i=0}^{l-1} \frac{(yl^2)^i}{(i!)^2} \le 2^{l-1}(y l^3) \max_{0 \le i \le l-1} \frac{(yl^2)^i}{(i!)^2}. \]
By quotient test we find that the maximal value is attained at $i_0=[\sqrt{y}l]$. If $i_0=0$ or $1$, then the equality (\ref{3:bmp}) can be easily verified. Now suppose $i_0 \ge 2$. Then $i_0 >\sqrt{y} r-1 \ge \sqrt{y} r/2$. Using the Stirling's bound on $n!$, given by
\begin{eqnarray} \label{3:stir} n! \ge \sqrt{2 \pi n} (n/e)^n, \end{eqnarray}
we obtain
\[\left|b_{\MP}^{(l)}\right|< 2^{l-1}(y l^3) \frac{(yl^2)^{i_0}}{4 \pi \left(\sqrt{y}l/2e\right)^{2i_0}}=\frac{l^32^l }{8 \pi}(4e^2)^{i_0}y \le \frac{l^32^l }{8 \pi}(4e^2)^{l}y. \]
This completes the proof of Lemma \ref{3:lem1}. \quad $\square$

To prove Theorem \ref{1:thm}, following the method of \cite{ba1}, we need a lemma from probability theory, which is discussed in details in \cite[Ch. XVI-3]{fel} (or see \cite[Lemma 3.1]{ba1}):

\begin{lemma} \label{3:prob} Let $F$ be a probability distribution with vanishing expectation and characteristic function $\phi$. Suppose that $F-G$ vanishes at $\pm \infty$ and that $G$ has a derivative $g$ such that $|g| \le m$. Finally, suppose that $g$ has a continuously differentiable Fourier transform $\gamma$ such that $\gamma(0)=1$ and $\gamma'(0)=0$. Then, for all $z$ and $T>0$ we have
\[|F(z)-G(z)| \le \frac{1}{\pi} \int_{-T}^T \left|\frac{\phi(t)-\gamma(t)}{t}\right|\, \ud t+\frac{24 m}{\pi T}. \]
\end{lemma}

\subsection{Proof of Theorem \ref{1:thm}} Using notation from Section \ref{sec-moment}, for each $s \in \Omega_p$, let $\lambda_1(s), \ldots, \lambda_p(s)$ be the eigenvalues of $\G$. The characteristic function we consider is
\[\phi_{\C}(t):=\frac{1}{p} \sum_{k=1}^p \E\left(\exp\left(it (\lambda_k(s)-1)\right),\Omega_p\right).\]
For the Marchenko-Pastur random variable $\vec{x}$ we consider
\[\gamma(t):=\E\left(\exp(it (\vec{x}-1) )\right). \]
Define for each $l$
\begin{eqnarray*} B_l=\frac{1}{p}\sum_{k=1}^p \E \left((\lambda_k(s)-1)^l, \Omega_p\right). \end{eqnarray*}
Expanding the $l$-th power we find that
\begin{eqnarray} \label{3:bl} B_l=\sum_{t=0}^l (-1)^{l-t} \binom{l}{t} \E \left(A_t(s), \Omega_p\right), \end{eqnarray}
where estimates on $\E \left(A_t(s), \Omega_p\right)$ is provided by Theorem \ref{2:moment}. Using the inequality
\[\left|\exp(it)-\sum_{l=0}^{r-1}\frac{(it)^l}{l!}\right| \le \frac{|t|^r}{r!}, \]
and choosing $r \ge 4$ to be even, we find that
\begin{eqnarray} \label{3:phi} \left|\phi_{\C}(t)-\sum_{l=0}^{r-1}\frac{(it)^l B_l}{l!}\right| \le \frac{t^r B_r}{r!},\end{eqnarray}
and
\begin{eqnarray} \label{3:gamma} \left|\gamma(t)-\sum_{l=0}^{r-1}\frac{(it)^l b_{\MP}^{(l)}}{l!}\right| \le \frac{t^r \, b_{\MP}^{(r)} }{r!} .\end{eqnarray}
We note that $B_l=b_{\MP}^{(l)}$ for $l=0,1$. For $l \ge 2$, using the expression (\ref{3:bl}) and Theorem \ref{2:moment}, given that $d^{\bot} \ge 5$, we find
\begin{eqnarray} \label{3:bll} \left|B_l-b_{\MP}^{(l)}\right| \le \sum_{t=2}^l \binom{l}{t} \frac{5 \, t^{t+1}}{n} <\frac{15 \, l^{l+1}}{n}. \end{eqnarray}
In writing
\[|\phi_{\C}(t)-\gamma(t)| \le \left|\phi_{\C}(t)-\sum_{l=0}^{r-1}\frac{(it)^l B_l}{l!}\right|+\left|\gamma(t)-\sum_{l=0}^{r-1}\frac{(it)^l b_{\MP}^{(l)}}{l!}\right|+\left|\sum_{l=0}^{r-1}\frac{(it)^l
\left(B_l-b_{\MP}^{(l)}\right)}{l!}\right|, \]
applying Lemma \ref{3:prob} and using the above estimates from (\ref{3:phi})(\ref{3:gamma})(\ref{3:bll}) and Lemma \ref{3:lem1}, we collect terms together and finally obtain
\begin{eqnarray} \label{3:mcest} \left|M_{\C}(z+1)-M_{\MP}(z+1)\right| \le \frac{r^2(8e^2T)^r y}{2 \pi^2 (r!)}+\frac{60 \, r (Tr)^r}{\pi n (r!)}+ \frac{24}{\pi^2 \sqrt{y}(1-y) T}. \end{eqnarray}

Finally, taking $r$ to be a positive even integer of size
\[ r \approx \frac{\log n}{\log \log n}, \mbox{ and } \,\, T =\frac{r}{16e^3},\]
and using the Stirling's bound (\ref{3:stir}), when $n$ (and consequently $r$) is sufficiently large, it is easy to see that the first two terms on the right side of (\ref{3:mcest}) can be both bounded by $\frac{\log \log n}{\log n}$, while the third term is
\[\frac{24 \cdot 16 \cdot e^3}{\pi^2 \sqrt{y}(1-y) r} \approx \frac{782 \cdot \log \log n }{\sqrt{y}(1-y) \cdot \log n}. \]
Combining these terms completes the proof of Theorem \ref{1:thm}. \quad $\square$

\section{The analysis of $W_{\gamma}$} \label{sec-combin}

Let $\gamma: [0,l_{\gamma}] \to [1,p]$ be a closed path with $V_{\gamma}=\gamma([0,l_{\gamma}])=\{z_a: 1 \le a \le v_{\gamma}\}, v_{\gamma}=|V_{\gamma}|$ and $I_a=\gamma^{-1}(z_a)$. Denote by $W_{\gamma}$ the number of solutions $(t_0,t_1,\ldots,t_{l_{\gamma}-1})$ such that $1 \le t_0,t_1,\ldots,t_{l_{\gamma}-1} \le n$ and
\[\sum_{u \in I_a} \left(\vec{h}_{t_u}-\vec{h}_{t_{u-1}} \right)=\vec{0}, \quad \forall \, \, 1 \le a \le v_{\gamma},\]
here $\vec{h}_i$ denotes the $i$-th row of the matrix $H$, whose rows are all distinct by assumption, and the indices shall be considered modulo $l_{\gamma}$, i.e., $t_{-1}=t_{l_{\gamma}-1}$. The purpose of this section is to study $W_{\gamma}$, which is crucial in the proof of Theorem \ref{2:moment}.

\begin{defin}
The closed path $\gamma$ is called ``reduced'' if $v_{\gamma}=l_{\gamma}=1$, or if $v_{\gamma} \ge 2$ and the following two conditions are satisfied:
\begin{itemize}
\item[(i).] each $|I_a| \ge 2$, hence $l =\sum_a |I_a| \ge 2 v \ge 4$;

\item[(ii).] each $I_a$ does not contain consecutive indices, that is, $\gamma(u) \ne \gamma(u+1) \,, \forall u$.

\end{itemize}
\end{defin}

We first study $W_{\gamma}$ when $\gamma$ is reduced.

\subsection{Study of $W_{\gamma}$ for $\gamma$ reduced}

Let $\gamma$ be a reduced closed path with $l=l_{\gamma} \ge 1$ and $v=v_{\gamma} \ge 1$. If $v_{\gamma}=l_{\gamma}=1$, then trivially we have
\[W_{\gamma}=n. \]
Now suppose that $v_{\gamma} \ge 2$. For each $I_a$, define $I_a':=I_a-\{1\}=\left\{u-1 \pmod{l_{\gamma}-1}: u \in I_a\right\}$. For any $1 \le a \le v_{\gamma}$, the equation corresponding to $I_a$ is
\begin{eqnarray} \label{3:matrix} \sum_{u \in I_a} \vec{h}_{t_u} -\sum_{u \in I_a'} \vec{h}_{t_u}=\vec{0}. \end{eqnarray}
We shall write down the equations (\ref{3:matrix}) for $1 \le a \le v_{\gamma}$ as a matrix with respect to the variables $\vec{h}_{t_0},\vec{h}_{t_1}, \ldots,\vec{h}_{t_{l-1}}$, given in the same ordered.

Since $\cup_a I_a$ is a partition of $[0,l-1]$, and each $I_a$ does not contain consecutive elements, there are distinct indices, which we may say $1$ and $v$, such that $0 \in I_1$ and $1 \in I_v$. Hence $k-1 \in I_0'$, and the row vector corresponding to the equation of $I_1$ with respect to $\vec{h}_{t_0},\vec{h}_{t_1}, \ldots,\vec{h}_{t_{l-1}}$ is of shape
\[[1,*, \cdots, *, -1]. \]
Now let $u_2$ be the smallest index in the set $\cup_{2 \le a \le v-1}\left(I_a \cup I_a'\right)$. We must have $u_2 \ge 1$, and $u_2 \in I_a'$ for some $2 \le a \le v-1$, because if otherwise, then $u_2=0$, which contradicts the fact that $0 \in I_1$ and $1 \in I_v$. We may reorder the indices and say $u_2 \in I_2'$. Hence $u_2+1 \in I_2$, and the row vector corresponding to the equation of $I_2$ with respect to $\vec{h}_{t_0},\vec{h}_{t_1}, \ldots,\vec{h}_{t_{l-1}}$ is of shape
\[[0 \cdots 0,-1,1,*, \cdots, *, 0], \]
where the first non-zero entry ``$-1$'' appears at the $u_2$-th column.

Now let $u_3$ be the smallest index in the set $\cup_{3 \le a \le v-1}\left(I_a \cup I_a'\right)$. Similarly we must have $u_3 \ge u_2+1$, and $u_3 \in I_a'$ for some $3 \le a \le v-1$. We reorder the indices and say $u_3 \in I_3'$. Then $u_3+1 \in I_3$, and the row vector corresponding to the equation of $I_3$ with respect to $\vec{h}_{t_0},\vec{h}_{t_1}, \ldots,\vec{h}_{t_{l-1}}$ is of shape
\[[0 \cdots 0,0 \cdots 0,-1,1,*, \cdots, *, 0], \]
where the first non-zero entry ``$-1$'' appears at the $u_3$-th column.

We can continue this process up to $a=v-1$ because each row contains at least two non-zero entries. Clearly the row vectors corresponding to the equations $I_a$ for $1 \le a \le v-1$ form an upper triangular matrix with rank $v-1$. So the number of free variables is $l-v+1$. This proves that $W_{\gamma} \le n^{l-v+1}$. Actually we shall do much better.

Since $l \ge 2 v$, and each row vector corresponding to $I_a, 1 \le a \le v-1$ with respect to $\vec{h}_{t_0},\vec{h}_{t_1}, \ldots,\vec{h}_{t_{l-1}}$ contains at least two $1$'s, we may find $l-v$ free variables, say they are $t_{v},\ldots,t_{l-1}$ after reordering the indices, so that for any given values of $t_{v}, \ldots, t_{l-1}$ from $1$ to $n$, solving the equations (\ref{3:matrix}) becomes looking for $1 \le t_0,\ldots,t_{v-1} \le n$ such that
\begin{eqnarray*} \vec{h}_{t_i} &=& \vec{v}_i, \quad \forall \, 2 \le i \le v-1,\\
\vec{h}_{t_0}+\vec{h}_{t_1}&=&\vec{v}_1, \end{eqnarray*}
where the vectors $\vec{v}_i$ are linear combinations of the rows of $H$, depending only on $t_v,\ldots,t_{l-1}$. Clearly the number of solutions for $t_i, 2 \le i \le v-1$ is at most one. One only needs to consider $t_0,t_1$.

If $\vec{v}_1=\vec{0}$, this enforces a new relation on $t_v,\ldots,t_{l-1}$ which were free before, hence the number of such $(t_v,\ldots,t_{l-1})$'s with $\vec{v}_1=\vec{0}$ is at most $n^{l-v-1}$. On the other hand, for each given $t_0$, there is at most one value $t_1$ such that $\vec{h}_{t_0}+\vec{h}_{t_1}=\vec{0}$. Hence the total number of solutions of $t_i$'s for this case is at most $n^{l-v}$. Let us define
\[A_{\vec{v}}=|\left\{(t_0,t_1): 1 \le t_0,t_1 \le n, \mbox{ and } \vec{h}_{t_0}+\vec{h}_{t_1}=\vec{v} \right\}|. \]
We have just proved that
\begin{eqnarray} \label{3:wr} W_{\gamma} \le n^{l-v} \left(1+\sup_{\vec{v} \ne \vec{0}}A_{\vec{v}}\right). \end{eqnarray}
Now for a fixed $\vec{v} \ne \vec{0}$, note that if $t_0=t_1$, the equation $2 \vec{h}_{t}=\vec{v}$ has at most one solution for $1 \le t \le n$. So we have
\begin{eqnarray} \label{3:av} A_{\vec{v}} \le 1 + 2 B_{\vec{v}},  \end{eqnarray}
where $B_{\vec{v}}$ is the cardinality of the set
\[\B_{\vec{v}}=\left\{(t_0,t_1): 1 \le t_0<t_1 \le n, \mbox{ and } \vec{h}_{t_0}+\vec{h}_{t_1}=\vec{v} \right\}. \]
If $B_{\vec{v}} \ge 2$, then for any distinct elements $(t_0,t_1),(t_0',t_1') \in \B_{\vec{v}}$, we conclude that $t_0,t_1,t_0',t_1'$ are all distinct and
\[\vec{h}_{t_0}+\vec{h}_{t_1}-\vec{h}_{t_0'}-\vec{h}_{t_1'}=\vec{0}. \]
This gives a weight $4$ codeword in $\C^{\bot}$ with entries $1,1,-1,-1$ at the $t_0,t_1,t_0'$ and $t_1'$-th places respectively. From it we may multiply elements of $\F(q)-\{0\}$ to get new weight $4$ codewords. Now suppose that $A$ is the number of weight $4$ codewords of $\C^{\bot}$. The above argument shows that
\[A \ge (q-1) \binom{B_{\vec{v}}}{2}= \frac{q-1}{2} B_{\vec{v}} \left(B_{\vec{v}}-1\right). \]
Hence we have
\[B_{\vec{v}} \le \sqrt{\frac{2A}{q-1}+\frac{1}{4}}+\frac{1}{2}\,\,. \]
In relation to (\ref{3:av}) and (\ref{3:wr}) we conclude that if $v_{\gamma} \ge 2$,
\begin{eqnarray} \label{2:estwr} W_{\gamma} \le C_A \, n^{l_{\gamma}-v_{\gamma}},\end{eqnarray}
where
\begin{eqnarray} \label{2:ca}
C_A=3+2 \sqrt{\frac{2A}{q-1}+\frac{1}{4}}\,.
\end{eqnarray}

\subsection{An example}

To illuminate the combinatorial nature of solving $W_{\gamma}$ in general, it may be useful to consider an example first.

Let $l_{\gamma}=9$, and $\gamma$ define the partition
\[\{0,1,\ldots,8\}= \{0,1,2,7\} \cup \{3,5,8\} \cup \{4\} \cup \{6\}. \]
So $v_{\gamma}=4$. Then $W_{\gamma}$ is the number of solutions $(t_0,t_1,\ldots,t_8)$ such that $1 \le t_0,t_1,\ldots,t_8 \le n$ and the following four equations hold simultaneously:
\begin{eqnarray} \label{2:exm1} \vec{h}_{t_0}+\vec{h}_{t_1}+\vec{h}_{t_2}+\vec{h}_{t_7}&=&
\vec{h}_{t_8}+\vec{h}_{t_0}+\vec{h}_{t_1}+\vec{h}_{t_6}\\
\label{2:exm2} \vec{h}_{t_3}+\vec{h}_{t_5}+\vec{h}_{t_8}&=&
\vec{h}_{t_2}+\vec{h}_{t_4}+\vec{h}_{t_7}\\
\label{2:exm3} \vec{h}_{t_4}&=&
\vec{h}_{t_3}\\
\label{2:exm4} \vec{h}_{t_6}&=&
\vec{h}_{t_5}\end{eqnarray}
Clearly one equation is redundant: we can always remove one and keep the rest.

Consider (\ref{2:exm1}), we find that $\vec{h}_{t_0},\vec{h}_{t_1}$ can be canceled out on both sides. Hence $t_0$ and $t_1$ are free and can be removed, and (\ref{2:exm1}) becomes
\begin{eqnarray} \label{2:exm11} \vec{h}_{t_2}+\vec{h}_{t_7}&=&
\vec{h}_{t_8}+\vec{h}_{t_6}
\end{eqnarray}

Consider (\ref{2:exm3}). Since the rows of $H$ are all distinct, this implies that $t_3=t_4$, and under this restriction, $\vec{h}_{t_3}$ and $\vec{h}_{t_4}$ are also canceled out on both sides of (\ref{2:exm2}). Then $t_3=t_4$ is also a free variable and can be removed.

Consider (\ref{2:exm4}). Clearly we have $t_5=t_6$, but this is not a free variable: replacing $t_5$ by $t_6$, we find that $W_{\gamma}=n^3 W_{\gamma'}$, where $W_{\gamma'}$ is the number of solutions $(t_2,t_6,t_7,t_8)$ such that $1 \le t_2,t_6,t_7,t_8 \le n$ and the equation (\ref{2:exm11}) is satisfied.

The $\gamma'$ can be reinterpreted as a closed path. It is a reduced path with $l_{\gamma'}=4,v_{\gamma'}=2$, hence the quantity $W_{\gamma'}$ can be estimated by (\ref{2:estwr}), so we conclude that
\[W_{\gamma} \le n^3 \, C_A \, n^{l_{\gamma'}-v_{\gamma'}}=C_A \, n^{5}. \]

\subsection{Study of $W_{\gamma}$ in general}

As illustrated by the previous example, we shall isolate variables from the equations related to $W_{\gamma}$, and removing these variables would result in a new but simpler closed path $\gamma'$, and three different situations may arise and need to be examined carefully.

We use some notation. For a closed path $\gamma: [0,l_{\gamma}] \to [1,p]$, the terms $V_{\gamma}, v_{\gamma}$ and $I_a$'s are as before. $\gamma$ yields a loop $t_0, t_1, \cdots, t_{u-1},t_u,t_{u+1}, \cdots, t_{l-2},t_{l-1},t_0 $, according to which we say that $t_{u-1}$ and $t_{u}$ are consecutive in $\gamma$, and $t_u':=t_{u-1}$ is the left neighbor of $t_u$ (as usual $t_{l-1}$ is the left neighbor of $t_0$). If we remove $t_u$ from $\gamma$, then in the resulting $\gamma'$, the loop is $t_0,\ldots,t_{u-1},t_{u+1},\ldots,t_{l_{\gamma}-1},t_0$, hence $l_{\gamma'}=l-1$, and the left neighbor of $t_{u+1}$ becomes $t_{u-1}$, but all other relations in terms of ``left neighbors'' stay the same.

\subsubsection{Case 1. Removing consecutive elements}

Suppose that there are consecutive elements in $I_a$ for some $a$, say, for example $u,u+1 \in I_a$. The equation with respect to $I_a$ is
\[\cdots+\vec{h}_{t_u}+\vec{h}_{t_{u+1}}+\cdots=
\cdots+\vec{h}_{t_{u-1}}+\vec{h}_{t_{u}}+\cdots.\]
Clearly $\vec{h}_{t_u}$ can be canceled out on both sides of the equation, and it does not appear in any other equations with respect to $I_b$, $b \ne a$. Let $\gamma'$ be the closed path by removing $t_u$, then $t_{u-1}$ becomes the left neighbor of $t_{u+1}$ in $\gamma'$ and all other relations in terms of ``neighbors'' remain the same. Hence we have
\[ \mbox{ Case 1}: \qquad l_{\gamma'}=l-1, \,\,\, v_{\gamma'}=v_{\gamma},\,\,\, W_{\gamma}=nW_{\gamma'}. \]
In $W_{\gamma'}$, we may rename the variables so that $\gamma':[0,l_{\gamma'}] \to [1,p]$ is a closed path with variables $t_0, \ldots,t_{l_{\gamma'-1}}$.

\subsubsection{Case 2. Removing ``leaves''} For a closed path $\gamma$, the vertex $u \in I_a$ is called a ``leaf'' if $I_a=\{u\}$ and $\gamma(u-1)=\gamma(u+1) \ne \gamma(u)$. Hence $u-1,u+1 \in I_b$ for some $b \ne a$. The equation with respect to $I_a$ is
\begin{eqnarray} \label{2:arg1} \vec{h}_{t_u}=\vec{h}_{t_{u-1}} \Longrightarrow t_u=t_{u-1}. \end{eqnarray}
The equation with respect to $I_b$ is
\begin{eqnarray} \label{2:arg2} \cdots+\vec{h}_{t_{u-1}}+\vec{h}_{t_{u+1}}+ \cdots= \cdots+\vec{h}_{t_{u-2}}+\vec{h}_{t_{u}}+ \cdots. \end{eqnarray}
Assuming (\ref{2:arg1}), then $\vec{h}_{t_u}$ and $\vec{h}_{t_{u-1}}$ can be canceled out trivially on both sides of (\ref{2:arg2}). Hence we have solved that $t_u=t_{u-1}$, which can be removed from the variables. Let $\gamma'$ be the resulting closed path. Removing both $t_u,t_{u-1}$ from (\ref{2:arg1}), it is clear that in $\gamma'$, $t_{u-2}$ becomes the left neighbor of $t_{u+1}$ and all other relations in terms of ``neighbors'' remain the same. We have
\[ \mbox{ Case 2}: \qquad l_{\gamma'}=l-2, \,\,\, v_{\gamma'}=v_{\gamma}-1,\,\,\, W_{\gamma}=nW_{\gamma'}. \]

\subsubsection{Case 3. Removing ``transition'' vertices} For a closed path $\gamma$, the vertex $u \in I_a$ is called a ``transition'' vertex if $I_a=\{u\}$ and $\gamma(u-1), \gamma(u),\gamma(u+1)$ are all distinct. Say $u-1 \in I_b$ and $u+1 \in I_c$, where $a,b,c$ are all distinct. The equation with respect to $I_a$ is still
\begin{eqnarray} \label{2:arg31} \vec{h}_{t_u}=\vec{h}_{t_{u-1}} \Longrightarrow t_u=t_{u-1}. \end{eqnarray}
The equations with respect to $I_b,I_c$ are
\begin{eqnarray} \label{2:arg32} \cdots+\vec{h}_{t_{u-1}}+\cdots &=& \cdots+\vec{h}_{t_{u-2}}+ \cdots \\
 \label{2:arg33} \cdots+\vec{h}_{t_{u+1}}+\cdots &=& \cdots+\vec{h}_{t_{u}}+ \cdots  \end{eqnarray}
Assuming (\ref{2:arg31}), that is, replacing $t_u$ by $t_{u-1}$, then (\ref{2:arg32}) stays the same but (\ref{2:arg33}) becomes
\begin{eqnarray*} \label{2:arg34} \cdots+\vec{h}_{t_{u+1}}+\cdots &=& \cdots+\vec{h}_{t_{u-1}}+ \cdots  \end{eqnarray*}
which means that by removing $t_u$, in the resulting $\gamma'$, $t_{u-1}$ becomes the left neighbor of $t_{u+1}$ and all the other relations in terms of ``neighbors'' remain the same. So we have
\[ \mbox{ Case 3}: \qquad l_{\gamma'}=l-1, \,\,\, v_{\gamma'}=v_{\gamma}-1,\,\,\, W_{\gamma}=W_{\gamma'}. \]

\subsection{Conclusion on $W_{\gamma}$}

In conclusion, suppose that altogether we perform $u,v$, and $w (\ge 0)$ times of Case 1, Case 2 and Case 3 reductions respectively on $\gamma$, maybe in different orders and combinations, to finally arrive at, after reordering the variables, a closed path $\gamma':[0,l_{\gamma'}] \to [1,p]$ with $l_{\gamma'},v_{\gamma'} \ge 1$, on which we could not do any of the reductions as described above. Then by definition $\gamma'$ is a reduced path, and we also have
\begin{eqnarray} \label{3:genwr} l_{\gamma'}=l_{\gamma}-u-2v-w , \,\,\, v_{\gamma'}=v_{\gamma}-v-w ,\,\,\, W_{\gamma}=n^{u+v}W_{\gamma'}\,\,. \end{eqnarray}
There are two cases:

\noindent {\bf Case 1.} If $v_{\gamma'}=l_{\gamma'}=1$, then $W_{\gamma'}=n$. Hence in this case $W_{\gamma}=n^{l_{\gamma}-v_{\gamma}+1}$.

\noindent {\bf Case 2.} If $v_{\gamma'} \ge 2$, then $W_{\gamma'} \le C_A \, n^{l_{\gamma'}-v_{\gamma'}}$ by (\ref{2:estwr}). We have in this case $W_{\gamma} \le C_A \, n^{l_{\gamma}-v_{\gamma}}$.

Denote by $\Gamma$ the set of all the $\gamma$'s that can be reduced to {\bf Case 1}. We conclude that
\[\left\{\begin{array}{lll} W_{\gamma}= n^{l_{\gamma}-v_{\gamma}+1} &: & \mbox{ if } \gamma \in \Gamma; \\
0 \le W_{\gamma} \le C_A \, n^{l_{\gamma}-v_{\gamma}} &: & \mbox{ if } \gamma \notin \Gamma, \end{array}
\right. \]

\subsection{Combinatorial structure of $\Gamma$}

Finally we need to prove the identity
\begin{eqnarray} \label{4:graph} \sum_{\substack{\gamma \in \Gamma \subset \Pi_{l,p}/\Sigma_p\\
v_{\gamma}=v}} 1=\frac{1}{v} \binom{l}{v-1}\binom{l-1}{v-1}. \end{eqnarray}
The theory of random matrices has been extensively studied (see \cite{and,meh}), and the above identity might be a well-known fact. Actually the left hand side appears naturally in the standard proof of the Marchenko-Pastur law for random matrices. Since we can not find a reference, we may sketch a proof here.

Let $X=(\vec{x}_{ij}) \in \R^{p \times n}$ be a random matrix where $\vec{x}_{ij}$'s are i.i.d, $\E(\vec{x}_{ij}) = 0,\E(\vec{x}^2_{ij}) = 1$ and $p<n$. Define
\[S=\frac{1}{n}XX^T.\]
Then
\[\frac{1}{p}\,\, \E\left(\tr(S^l)\right)=\frac{1}{pn^l}\sum_{\gamma,\tau}\E\left(
\vec{x}_{\gamma(0)t_0}\vec{x}_{\gamma(1)t_0}\vec{x}_{\gamma(1)t_1}\vec{x}_{\gamma(2)t_1} \ldots \vec{x}_{\gamma(l-1)t_{l-1}}\vec{x}_{\gamma(0)t_{l-1}}\right)=
\frac{1}{pn^l}\sum_{\gamma,\tau}\E\left(\gamma,\tau\right), \]
where the sum is over all maps $\gamma \in \Pi_{l,p}$ and all $\tau: =\{t_i\}_{i=0}^{l-1} \in [1,n]^l$. Now this corresponds to a directed loop on a bipartite graph from the vertex set $\{\gamma(0),\ldots,\gamma(l-1)\}$ to the vertex set $\{t_0,\ldots,t_{l-1}\}$ with $2l$ steps. As the standard proof goes, each edge must appear at least twice, otherwise $\E(\gamma,\tau)=0$. Hence we have at most $l$ edges in the graph, and at most $l + 1$ vertices in the skeleton. The optimal situation, that is, graphs with exactly $l$ edges and $l+1$ vertices, or ``double trees'' will give the main contribution. Terms arising from other configuration of graphs are negligible and can be ignored. The standard result on counting such ``double trees'' is that, for each $1 \le v \le l$, the number of double tree shapes with $v$ vertices in $\gamma$ (i.e., $v_{\gamma}=v$) and $l-v+1$ vertices in $\tau$ is given by the right hand side of (\ref{4:graph}) (see \cite[page 20, Exercise 2.1.18]{and}). A little thought about properties of $\Gamma$ concludes that the left hand side of (\ref{4:graph}) also counts the total number of such double trees. The finishes the proof of the identity (\ref{4:graph}). \quad $\square$






\end{document}